% Plain TeX interface to graphicx package.
% David Carlisle

\input miniltx

% --- shims for using modern graphics-def files under plain TeX ---
% plain's \newdimen is \outer, which modern pdftex.def uses inside a
% macro argument; redefine without \outer (same expansion as plain.tex)
\def\newdimen{\alloc@1\dimen\dimendef\insc@unt}
% LaTeX kernel log-message helpers referenced by pdftex.def's \PackageInfo
\def\on@line{ on input line \the\inputlineno}

% -----------------------------------------------------------------

\def\Gin@driver{pdftex.def}
\input graphicx.sty

\resetatcatcode

\input eplain
\newfam\bffam
\textfont\bffam=\tenbf\scriptfont\bffam=\sevenbf
\scriptscriptfont\bffam=\fivebf
\def\bf{\fam\bffam\tenbf}
%=================================================================
% Blackboard Bold
%=================================================================
\def\bbbA{{\mathchoice 
         {\setbox0=\hbox{$\displaystyle{\hbox{\bf{A}}}$}\hbox{\raise
         0.15\ht0\hbox to0pt{\kern0.4\wd0\vrule height0.8\ht0\hss}\box0}}
         {\setbox0=\hbox{$\textstyle{\hbox{\bf{A}}}$}\hbox{\raise
         0.15\ht0\hbox to0pt{\kern0.4\wd0\vrule height0.8\ht0\hss}\box0}}
         {\setbox0=\hbox{$\scriptstyle{\hbox{\bf{A}}}$}\hbox{\raise
         0.15\ht0\hbox to0pt{\kern0.4\wd0\vrule height0.7\ht0\hss}\box0}}
         {\setbox0=\hbox{$\scriptscriptstyle{\hbox{\bf{A}}}$}\hbox{\raise
         0.15\ht0\hbox to0pt{\kern0.4\wd0\vrule height0.7\ht0\hss}\box0}}}}

\def\bbbC{{\mathchoice {\setbox0=\hbox{$\displaystyle{\hbox{\bf{C}}}$}\hbox{\raise
          0.15\ht0\hbox to0pt{\kern0.4\wd0\vrule height0.8\ht0\hss}\box0}}
          {\setbox0=\hbox{$\textstyle{\hbox{\bf{C}}}$}\hbox{\raise
          0.15\ht0\hbox to0pt{\kern0.4\wd0\vrule height0.8\ht0\hss}\box0}}
          {\setbox0=\hbox{$\scriptstyle{\hbox{\bf{C}}}$}\hbox{\raise
          0.15\ht0\hbox to0pt{\kern0.4\wd0\vrule height0.7\ht0\hss}\box0}}
          {\setbox0=\hbox{$\scriptscriptstyle{\hbox{\bf{C}}}$}\hbox{\raise
          0.15\ht0\hbox to0pt{\kern0.4\wd0\vrule height0.7\ht0\hss}\box0}}}}

\def\bbbG{{\mathchoice {\setbox0=\hbox{$\displaystyle{\hbox{\bf{G}}}$}\hbox{\raise
          0.15\ht0\hbox to0pt{\kern0.4\wd0\vrule height0.8\ht0\hss}\box0}}
          {\setbox0=\hbox{$\textstyle{\hbox{\bf{G}}}$}\hbox{\raise
          0.15\ht0\hbox to0pt{\kern0.4\wd0\vrule height0.8\ht0\hss}\box0}}
          {\setbox0=\hbox{$\scriptstyle{\hbox{\bf{G}}}$}\hbox{\raise
          0.15\ht0\hbox to0pt{\kern0.4\wd0\vrule height0.7\ht0\hss}\box0}}
          {\setbox0=\hbox{$\scriptscriptstyle{\hbox{\bf{G}}}$}\hbox{\raise
          0.15\ht0\hbox to0pt{\kern0.4\wd0\vrule height0.7\ht0\hss}\box0}}}}

\def\bbbJ{{\mathchoice {\setbox0=\hbox{$\displaystyle{\hbox{\bf{J}}}$}\hbox{\raise
          0.15\ht0\hbox to0pt{\kern0.4\wd0\vrule height0.8\ht0\hss}\box0}}
          {\setbox0=\hbox{$\textstyle{\hbox{\bf{J}}}$}\hbox{\raise
          0.15\ht0\hbox to0pt{\kern0.4\wd0\vrule height0.8\ht0\hss}\box0}}
          {\setbox0=\hbox{$\scriptstyle{\hbox{\bf{J}}}$}\hbox{\raise
          0.15\ht0\hbox to0pt{\kern0.4\wd0\vrule height0.7\ht0\hss}\box0}}
          {\setbox0=\hbox{$\scriptscriptstyle{\hbox{\bf{J}}}$}\hbox{\raise
          0.15\ht0\hbox to0pt{\kern0.4\wd0\vrule height0.7\ht0\hss}\box0}}}}

\def\bbbO{{\mathchoice {\setbox0=\hbox{$\displaystyle{\hbox{\bf{O}}}$}\hbox{\raise
          0.15\ht0\hbox to0pt{\kern0.4\wd0\vrule height0.8\ht0\hss}\box0}}
          {\setbox0=\hbox{$\textstyle{\hbox{\bf{O}}}$}\hbox{\raise
          0.15\ht0\hbox to0pt{\kern0.4\wd0\vrule height0.8\ht0\hss}\box0}}
          {\setbox0=\hbox{$\scriptstyle{\hbox{\bf{O}}}$}\hbox{\raise
          0.15\ht0\hbox to0pt{\kern0.4\wd0\vrule height0.7\ht0\hss}\box0}}
          {\setbox0=\hbox{$\scriptscriptstyle{\hbox{\bf{O}}}$}\hbox{\raise
          0.15\ht0\hbox to0pt{\kern0.4\wd0\vrule height0.7\ht0\hss}\box0}}}}

\def\bbbQ{{\mathchoice {\setbox0=\hbox{$\displaystyle{\hbox{\bf{Q}}}$}\hbox{\raise
          0.15\ht0\hbox to0pt{\kern0.4\wd0\vrule height0.8\ht0\hss}\box0}}
          {\setbox0=\hbox{$\textstyle{\hbox{\bf{Q}}}$}\hbox{\raise
          0.15\ht0\hbox to0pt{\kern0.4\wd0\vrule height0.8\ht0\hss}\box0}}
          {\setbox0=\hbox{$\scriptstyle{\hbox{\bf{Q}}}$}\hbox{\raise
          0.15\ht0\hbox to0pt{\kern0.4\wd0\vrule height0.7\ht0\hss}\box0}}
          {\setbox0=\hbox{$\scriptscriptstyle{\hbox{\bf{Q}}}$}\hbox{\raise
          0.15\ht0\hbox to0pt{\kern0.4\wd0\vrule height0.7\ht0\hss}\box0}}}}

\def\bbbS{{\mathchoice{\setbox0=\hbox{$\displaystyle{\hbox{\bf{S}}}$}
          \hbox{\raise0.5\ht0\hbox
          to0pt{\kern0.35\wd0\vrule height0.45\ht0\hss}\hbox
          to0pt{\kern0.55\wd0\vrule height0.5\ht0\hss}\box0}}
          {\setbox0=\hbox{$\textstyle {\hbox{\bf{S}}}$}\hbox{\raise0.5\ht0\hbox
          to0pt{\kern0.35\wd0\vrule height0.45\ht0\hss}\hbox
          to0pt{\kern0.55\wd0\vrule height0.5\ht0\hss}\box0}}
          {\setbox0=\hbox{$\scriptstyle{\hbox{\bf{S}}}$}\hbox{\raise0.5\ht0\hbox
          to0pt{\kern0.35\wd0\vrule height0.45\ht0\hss}\raise0.05\ht0\hbox
          to0pt{\kern0.5\wd0\vrule height0.45\ht0\hss}\box0}}
          {\setbox0=\hbox{$\scriptscriptstyle{\hbox{\bf{S}}}$}\hbox{\raise0.5\ht0\hbox
          to0pt{\kern0.4\wd0\vrule height0.45\ht0\hss}\raise0.05\ht0\hbox
          to0pt{\kern0.55\wd0\vrule height0.45\ht0\hss}\box0}}}}
\def\bbbT{{\mathchoice {\setbox0=\hbox{$\displaystyle\bf
          T$}\hbox{\hbox to0pt{\kern0.3\wd0\vrule height0.9\ht0\hss}\box0}}
          {\setbox0=\hbox{$\textstyle{\hbox{\bf{T}}}$}\hbox{\hbox
          to0pt{\kern0.3\wd0\vrule height0.9\ht0\hss}\box0}}
          {\setbox0=\hbox{$\scriptstyle{\hbox{\bf{T}}}$}\hbox{\hbox
          to0pt{\kern0.3\wd0\vrule height0.9\ht0\hss}\box0}}
          {\setbox0=\hbox{$\scriptscriptstyle{\hbox{\bf{T}}}$}\hbox{\hbox
          to0pt{\kern0.3\wd0\vrule height0.9\ht0\hss}\box0}}}}
\def\bbbU{{\mathchoice {\setbox0=\hbox{$\displaystyle{\hbox{\bf{U}}}$}\hbox{\raise
          0.15\ht0\hbox to0pt{\kern0.4\wd0\vrule height0.8\ht0\hss}\box0}}
          {\setbox0=\hbox{$\textstyle{\hbox{\bf{U}}}$}\hbox{\raise
          0.15\ht0\hbox to0pt{\kern0.4\wd0\vrule height0.8\ht0\hss}\box0}}
          {\setbox0=\hbox{$\scriptstyle{\hbox{\bf{U}}}$}\hbox{\raise
          0.15\ht0\hbox to0pt{\kern0.4\wd0\vrule height0.7\ht0\hss}\box0}}
          {\setbox0=\hbox{$\scriptscriptstyle{\hbox{\bf{U}}}$}\hbox{\raise
          0.15\ht0\hbox to0pt{\kern0.4\wd0\vrule height0.7\ht0\hss}\box0}}}}
\def\bbbV{{\mathchoice {\hbox{$\bf\textstyle V\kern-0.4em V$}}
          {\hbox{$\bf\textstyle V\kern-0.4em V$}}
          {\hbox{$\bf\scriptstyle V\kern-0.3em V$}}
          {\hbox{$\bf\scriptscriptstyle V\kern-0.2em V$}}}}
\def\bbbW{{\mathchoice {\hbox{$\bf\textstyle W\kern-0.4em W$}}
          {\hbox{$\bf\textstyle W\kern-0.4em W$}}
          {\hbox{$\bf\scriptstyle W\kern-0.3em W$}}
          {\hbox{$\bf\scriptscriptstyle W\kern-0.2em W$}}}}
\def\bbbX{{\mathchoice {\hbox{$\bf\textstyle X\kern-0.4em X$}}
          {\hbox{$\bf\textstyle X\kern-0.4em X$}}
          {\hbox{$\bf\scriptstyle X\kern-0.3em X$}}
          {\hbox{$\bf\scriptscriptstyle X\kern-0.2em X$}}}}
\def\bbbY{{\mathchoice {\hbox{$\bf\textstyle Y\kern-0.4em Y$}}
          {\hbox{$\bf\textstyle Y\kern-0.4em Y$}}
          {\hbox{$\bf\scriptstyle Y\kern-0.3em Y$}}
          {\hbox{$\bf\scriptscriptstyle Y\kern-0.2em Y$}}}}

\def\bbba{{\mathchoice {\setbox0=\hbox{$\displaystyle{\hbox{\bf{a}}}$}\hbox{\raise
          0.15\ht0\hbox to0pt{\kern0.4\wd0\vrule height0.8\ht0\hss}\box0}}
          {\setbox0=\hbox{$\textstyle{\hbox{\bf{a}}}$}\hbox{\raise
          0.15\ht0\hbox to0pt{\kern0.4\wd0\vrule height0.8\ht0\hss}\box0}}
          {\setbox0=\hbox{$\scriptstyle{\hbox{\bf{a}}}$}\hbox{\raise
          0.15\ht0\hbox to0pt{\kern0.4\wd0\vrule height0.7\ht0\hss}\box0}}
          {\setbox0=\hbox{$\scriptscriptstyle{\hbox{\bf{a}}}$}\hbox{\raise
          0.15\ht0\hbox to0pt{\kern0.4\wd0\vrule height0.7\ht0\hss}\box0}}}}
\def\bbbb{{\mathchoice {\setbox0=\hbox{$\displaystyle{\hbox{\bf{b}}}$}\hbox{\raise
          0.15\ht0\hbox to0pt{\kern0.4\wd0\vrule height0.8\ht0\hss}\box0}}
          {\setbox0=\hbox{$\textstyle{\hbox{\bf{b}}}$}\hbox{\raise
          0.15\ht0\hbox to0pt{\kern0.4\wd0\vrule height0.8\ht0\hss}\box0}}
          {\setbox0=\hbox{$\scriptstyle{\hbox{\bf{b}}}$}\hbox{\raise
          0.15\ht0\hbox to0pt{\kern0.4\wd0\vrule height0.7\ht0\hss}\box0}}
          {\setbox0=\hbox{$\scriptscriptstyle{\hbox{\bf{b}}}$}\hbox{\raise
          0.15\ht0\hbox to0pt{\kern0.4\wd0\vrule height0.7\ht0\hss}\box0}}}}
\def\bbbc{{\mathchoice {\setbox0=\hbox{$\displaystyle{\hbox{\bf{c}}}$}\hbox{\raise
          0.15\ht0\hbox to0pt{\kern0.4\wd0\vrule height0.8\ht0\hss}\box0}}
          {\setbox0=\hbox{$\textstyle{\hbox{\bf{c}}}$}\hbox{\raise
          0.15\ht0\hbox to0pt{\kern0.4\wd0\vrule height0.8\ht0\hss}\box0}}
          {\setbox0=\hbox{$\scriptstyle{\hbox{\bf{c}}}$}\hbox{\raise
          0.15\ht0\hbox to0pt{\kern0.4\wd0\vrule height0.7\ht0\hss}\box0}}
          {\setbox0=\hbox{$\scriptscriptstyle{\hbox{\bf{c}}}$}\hbox{\raise
          0.15\ht0\hbox to0pt{\kern0.4\wd0\vrule height0.7\ht0\hss}\box0}}}}
\def\bbbd{{\mathchoice {\setbox0=\hbox{$\displaystyle{\hbox{\bf{d}}}$}\hbox{\raise
          0.15\ht0\hbox to0pt{\kern0.4\wd0\vrule height0.8\ht0\hss}\box0}}
          {\setbox0=\hbox{$\textstyle{\hbox{\bf{d}}}$}\hbox{\raise
          0.15\ht0\hbox to0pt{\kern0.4\wd0\vrule height0.8\ht0\hss}\box0}}
          {\setbox0=\hbox{$\scriptstyle{\hbox{\bf{d}}}$}\hbox{\raise
          0.15\ht0\hbox to0pt{\kern0.4\wd0\vrule height0.7\ht0\hss}\box0}}
          {\setbox0=\hbox{$\scriptscriptstyle{\hbox{\bf{d}}}$}\hbox{\raise
          0.15\ht0\hbox to0pt{\kern0.4\wd0\vrule height0.7\ht0\hss}\box0}}}}
\def\bbbe{{\mathchoice {\setbox0=\hbox{$\displaystyle{\hbox{\bf{e}}}$}\hbox{\raise
          0.15\ht0\hbox to0pt{\kern0.4\wd0\vrule height0.8\ht0\hss}\box0}}
          {\setbox0=\hbox{$\textstyle{\hbox{\bf{e}}}$}\hbox{\raise
          0.15\ht0\hbox to0pt{\kern0.4\wd0\vrule height0.8\ht0\hss}\box0}}
          {\setbox0=\hbox{$\scriptstyle{\hbox{\bf{e}}}$}\hbox{\raise
          0.15\ht0\hbox to0pt{\kern0.4\wd0\vrule height0.7\ht0\hss}\box0}}
          {\setbox0=\hbox{$\scriptscriptstyle{\hbox{\bf{e}}}$}\hbox{\raise
          0.15\ht0\hbox to0pt{\kern0.4\wd0\vrule height0.7\ht0\hss}\box0}}}}
\def\bbbf{{\mathchoice {\setbox0=\hbox{$\displaystyle{\hbox{\bf{f}}}$}\hbox{\raise
          0.15\ht0\hbox to0pt{\kern0.4\wd0\vrule height0.8\ht0\hss}\box0}}
          {\setbox0=\hbox{$\textstyle{\hbox{\bf{f}}}$}\hbox{\raise
          0.15\ht0\hbox to0pt{\kern0.4\wd0\vrule height0.8\ht0\hss}\box0}}
          {\setbox0=\hbox{$\scriptstyle{\hbox{\bf{f}}}$}\hbox{\raise
          0.15\ht0\hbox to0pt{\kern0.4\wd0\vrule height0.7\ht0\hss}\box0}}
          {\setbox0=\hbox{$\scriptscriptstyle{\hbox{\bf{f}}}$}\hbox{\raise
          0.15\ht0\hbox to0pt{\kern0.4\wd0\vrule height0.7\ht0\hss}\box0}}}}
\def\bbbg{{\mathchoice {\setbox0=\hbox{$\displaystyle{\hbox{\bf{g}}}$}\hbox{\raise
          0.15\ht0\hbox to0pt{\kern0.4\wd0\vrule height0.8\ht0\hss}\box0}}
          {\setbox0=\hbox{$\textstyle{\hbox{\bf{g}}}$}\hbox{\raise
          0.15\ht0\hbox to0pt{\kern0.4\wd0\vrule height0.8\ht0\hss}\box0}}
          {\setbox0=\hbox{$\scriptstyle{\hbox{\bf{g}}}$}\hbox{\raise
          0.15\ht0\hbox to0pt{\kern0.4\wd0\vrule height0.7\ht0\hss}\box0}}
          {\setbox0=\hbox{$\scriptscriptstyle{\hbox{\bf{g}}}$}\hbox{\raise
          0.15\ht0\hbox to0pt{\kern0.4\wd0\vrule height0.7\ht0\hss}\box0}}}}
\def\bbbu{{\mathchoice {\setbox0=\hbox{$\displaystyle{\hbox{\bf{u}}}$}\hbox{\raise
          0.15\ht0\hbox to0pt{\kern0.4\wd0\vrule height0.8\ht0\hss}\box0}}
          {\setbox0=\hbox{$\textstyle{\hbox{\bf{u}}}$}\hbox{\raise
          0.15\ht0\hbox to0pt{\kern0.4\wd0\vrule height0.8\ht0\hss}\box0}}
          {\setbox0=\hbox{$\scriptstyle{\hbox{\bf{u}}}$}\hbox{\raise
          0.15\ht0\hbox to0pt{\kern0.4\wd0\vrule height0.7\ht0\hss}\box0}}
          {\setbox0=\hbox{$\scriptscriptstyle{\hbox{\bf{u}}}$}\hbox{\raise
          0.15\ht0\hbox to0pt{\kern0.4\wd0\vrule height0.7\ht0\hss}\box0}}}}
\def\bbby{{\mathchoice {\setbox0=\hbox{$\displaystyle{\hbox{\bf{y}}}$}\hbox{\raise
          0.15\ht0\hbox to0pt{\kern0.4\wd0\vrule height0.8\ht0\hss}\box0}}
          {\setbox0=\hbox{$\textstyle{\hbox{\bf{y}}}$}\hbox{\raise
          0.15\ht0\hbox to0pt{\kern0.4\wd0\vrule height0.8\ht0\hss}\box0}}
          {\setbox0=\hbox{$\scriptstyle{\hbox{\bf{y}}}$}\hbox{\raise
          0.15\ht0\hbox to0pt{\kern0.4\wd0\vrule height0.7\ht0\hss}\box0}}
          {\setbox0=\hbox{$\scriptscriptstyle{\hbox{\bf{y}}}$}\hbox{\raise
          0.15\ht0\hbox to0pt{\kern0.4\wd0\vrule height0.7\ht0\hss}\box0}}}}

\input miniltx

% twelvepoint.tex --- honest 12pt for plain TeX, with no \magnification.
%
% Purpose: \magnification's meaning has changed across pdftex
% generations (old pdftex scaled everything including the page;
% current pdftex scales content but not the page registers; some
% versions ignore it for output entirely). Any document using \mag
% therefore renders differently on different systems. This file
% produces 12pt output by loading fonts at 12pt and setting honest
% dimensions, so the SAME source gives the SAME output on every
% TeX from the 1990s to today, in DVI or PDF mode, locally or on
% arXiv.
%
% Usage: replace "\magnification=1200" with "\input twelvepoint".
% Do not use \magnification, \magstep, or "true" units afterward.

% --- text and math fonts at 12pt (design-size fonts where they exist) ---
\font\twelverm=cmr12         \font\twelvebf=cmbx12
\font\twelveit=cmti12        \font\twelvesl=cmsl12
\font\twelvett=cmtt12
\font\twelvei=cmmi12         \skewchar\twelvei='177
\font\twelvesy=cmsy10 at 12pt \skewchar\twelvesy='60
\font\twelveex=cmex10 at 12pt

% --- script (9pt) and scriptscript (7pt) sizes ---
\font\ninerm=cmr9            \font\ninebf=cmbx9
\font\ninei=cmmi9            \skewchar\ninei='177
\font\ninesy=cmsy9           \skewchar\ninesy='60
% (the 7pt fonts \sevenrm, \seveni, \sevensy, \sevenbf come from plain.tex)

% --- wire up the math families ---
\textfont0=\twelverm \scriptfont0=\ninerm \scriptscriptfont0=\sevenrm
\textfont1=\twelvei  \scriptfont1=\ninei  \scriptscriptfont1=\seveni
\textfont2=\twelvesy \scriptfont2=\ninesy \scriptscriptfont2=\sevensy
\textfont3=\twelveex \scriptfont3=\twelveex \scriptscriptfont3=\twelveex
\textfont\itfam=\twelveit
\textfont\slfam=\twelvesl
\textfont\bffam=\twelvebf \scriptfont\bffam=\ninebf
  \scriptscriptfont\bffam=\sevenbf
\textfont\ttfam=\twelvett

% --- size-switching macros ---
\def\rm{\fam0\twelverm}
\def\it{\fam\itfam\twelveit}
\def\sl{\fam\slfam\twelvesl}
\def\bf{\fam\bffam\twelvebf}
\def\tt{\fam\ttfam\twelvett}

% --- spacing scaled to 12pt ---
\normalbaselineskip=14.5pt
\setbox\strutbox=\hbox{\vrule height10.2pt depth4.2pt width0pt}
\normalbaselines
\jot=3.6pt

% --- honest page dimensions (no "true" units anywhere) ---
\hsize=6.5in
\vsize=8.9in

\rm
% end of twelvepoint.tex

\footline{\hfil\folio\hfil\the\month / \the\day / \the\year}
\def\qed{\line{\hfill q.e.d.}}
\parskip=8pt
\overfullrule=0pt
\newcount\secno  \newcount\itemno  \newcount\dnumno
\secno=0
\def\newsec#1{
                \advance\secno by 1
                \dnumno=0
                \itemno=0
                \bigskip 
                \noindent{\bf \the\secno.  #1\par}
                }

   \def\bigo#1{O\Big(#1\Big)}   
        
\newcount\prevsecno
\def\prevsec{\prevsecno=\secno\advance\prevsecno by -1 \the\prevsecno}
\def\prevprevsec{\prevsecno=\secno\advance\prevsecno by -2 \the\prevsecno}
\def\num{{\global\advance\itemno by 1}{\bf \the\secno.\the\itemno}}
\def\dnum{{\global\advance\dnumno by 1 }\eqno{(\the\secno.\the\dnumno})}
\centerline{\bf Some Integer Sequences Arising From Game Theory}
\medskip
\centerline{\bf by}
\medskip
\centerline{\bf Steven E.  Landsburg}
\centerline{\bf University of Rochester}

\newcount\fnno\fnno=0
\def\fn#1{\advance\fnno by 1 \footnote{$^\the\fnno$}{#1}}

\def\entry#1{\noindent {\bf #1 \num.}}
\secno=0

Alice and Bob are ice cream vendors.  Each announces a location; payoffs depend only on those announcements.   Specify the payoff functions and you have a good homework problem for Week One of a first course in game theory.

Actually, you have three good homework problems (or one good three-part homework problem):  What if Alice goes first?  What if Bob goes first?  What if they act simultaneously?  (Of course if the payoffs are symmetric, the three problems collapse to two.)

While I was constructing just such a homework set, it recently occurred to me that I could make the problem more interesting by having Alice and Bob each announce both a location and a price, with payoffs depending on all four announcements.   Now find all equilibria among all possible protocols.  (A protocol is something like ``First Alice announces a price, then they simultaneously announce locations, then Bob announces a price''.  Later in the paper, we will call these {\sl Hotelling protocols\/}, to distinguish them from the {\sl simple protocols\/} introduced below).

I don't always attempt to solve problems myself before I assign them, which doesn't always turn out well.  Fortunately, in this case I made a  partial exception to that rule by computing that there are 51 possible protocols (reducing to 27 if the payoffs are symmetric.)  This sufficed to deter me from assigning this problem.  It did not, however, deter me from wondering what happens to the number of protocols when Alice and Bob each have $n$ moves.  Nor did it deter me from spending far more time on this question than I can possibly justify.   This paper is the result.

The solution turns out to depend on solving what looks like a simpler problem:  Consider a game where all moves are of the same type; there is no such thing as a ``price move'' or a ``location move''.   At each stage of the game, either Alice makes one move or Bob makes one move or they each move simultaneously.  A {\it simple protocol\/} dictates this sequence in advance.    If Alice has $i$ moves and Bob has $j$ moves, how many simple protocols are there?    (For readers who want to check their understanding, the answer in the case $i=j=2$ is 9.)   
 
It turns out that counting simple protocols is a prerequisite for counting Hotelling protocols, so most of the paper (titled ``Part One'') deals with the simple case.
 In Part Two we return to the ice cream vendors by estimating the growth rate of $H(n)$, the number of Hotelling protocols when Alice and Bob have $n$ moves each.   Spoiler alert:  It grows very fast.  In fact, we will show that $H(n)$ grows like $(n!)^2(5.305)^n$.

{\bf A Remark on Methods.} There is a substantial literature on analytic methods for estimating power series coefficients (e.g. [PW] and [FS]), some of which is likely to be applicable here.   I made some brief attempts along these lines,  ran into difficulties of the sort documented in [MO], realized I'd stopped having fun, and decided to take a different tack.  Most of those techniques start with a fixed ratio $i/j$, whereas I had a strong preference  for fixing the difference $i-j$, and I chose to indulge that preference.  This preference was guided partly by personal aesthetics and partly by a desire to avoid number-theoretic complications such as having to worry about the common divisors of $i$ and $j$. Therefore there will be no analytic combinatorics in Part One, and for that matter, with just one exception, no analysis of any kind.  (The one exception is an elementary invocation of Cauchy's Theorem  in Section 3, which I suspect I could have avoided.)  Part Two requires (I suspect  less avoidably) an appeal to one key analytic result in the form of Temme's approximation for Stirling numbers of the second kind ([T]).  I do not know whether a more mainstream approach would yield better results than those obtained here, for any of various definitions of ``better''.  

{\bf A Remark on Notation.}   Notation established in numbered sections labeled ``Notation''   (e.g. the constant $C_1$ defined in Notation 5.1) is always in effect for the remainder of the paper.  Notation established in the midst of a proof or construction is always in local to that proof or construction.
 
\bigskip\bigskip

\centerline{\bf PART ONE:  SIMPLE PROTOCOLS}

\newsec{The Setup}

Alice and Bob play a game in which Alice moves $i$ times and Bob moves $j$ times.  At each stage of the game, either Alice makes one move, or Bob makes one move, or they each make one move simultaneously.  The possible sequences of play are enumerated by what we'll call $(i,j)$-strings:

\entry{Definition}  A {simple $(i,j)$-protocol\/} is a string of symbols where
\itemitem{a)}Each symbol is either $[A]$, $[B]$, or $[AB]$.
\itemitem{b)}No $[A]$ is followed by an $[A]$ and no $[B]$ is followed by a $[B]$.
\itemitem{c)}The number of $[A]$'s plus the number of $[AB]$'s is $i$
\itemitem{d)}The number of $[B]$'s plus the number of $[AB]$'s is $j$

\entry{Definition}
Let $t(i,j)$ be the number of simple $(i,j)$-protocols. We call these the $t$-numbers.

 The goal of Part One is to compute --- and to give simple estimates for --- the values of the $t(i,j)$ .  We can start with two elementary observations:
 
 \entry{Proposition}  a)  $t(i,j)=t(j,i)$.  
 
 b)  When $i>2j+1$, $t(i,j)=0$.
 
 {\bf Proof.}  a) is immediate.  b) is almost as immediate:  If Bob has $j$ moves, then Alice can have at most $j$ moves simultaneous with Bob and at most $j+1$ solo moves (because any two of her solo turns must be separated by a turn in which Bob moves), hence at most $2j+1$ moves altogether.

\qed

In the remainder of this section, we will derive an exact formula for $t(i,j)$, which unfortunately sheds little light on the growth rate and will not be used in the sequel.

\entry{Construction}
To construct an $(i,j)$-protocol with exactly $f$ occurrences of the symbol $[AB]$:

1)  Form a row of $2i+1$ boxes numbered $1,\ldots,2i+1$.

2)  Let $E=\{2,4,\ldots,2i\}$.  Define an $(f,r)$-subset of $E$ to be a subset $F\subset E$ of cardinality $f$, such that the number of runs of consecutive even integers in $E\backslash F$ is $r$. Choose $f$ and $r$ such that $0\le f \le j$ and $r\le f+1$.  

3)  Choose an  $(f,r)$ subset of $E$.  Label each of the corresponding boxes $[AB]$.

Note that the number of ways to complete step 3) is 

$$\pmatrix{f+1\cr r\cr}\pmatrix{i-f-1\cr r-1\cr}$$

4)  Label the remaining even-numbered boxes $[A]$.

5)  Call an (odd-numbered) box {\it forced\/} if it sits directly between two boxes labeled $[A]$.   Label all of the forced boxes $[B]$.  Note that there are exactly $i-f-r$ forced boxes.  (Proof: There are $i-f$ even boxes labeled $[A]$, arranged in $r$ runs.  Each run of length $\ell$ creates $\ell-1$ forced boxes.)

6)  Now there are $i+1-(i-f-r)=f+r+1$ unlabeled boxes.  There are $f+(i-f-r)=i-r$ boxes labeled either $[AB]$ or $[B]$.  Choose $j-(i-r)=r-i+j$ of the unlabeled boxes and label them $[B]$ (so that the number of $[AB]$ labels plus the number of $[B]$ labels is equal to $j$).  The number of ways to do this is 
$$\pmatrix{f+r+1\cr r-i+j\cr}$$

7)  Discard all unlabeled boxes. A few moments' thought should reveal that the sequence of remaining labels is an $(i,j)$-protocol, and that every $(i,j)$-protocol arises uniquely in this way. 
This suffices to compute $t(i,j)$, the number of $(i,j)$-protocols:

\entry{Theorem} For $i\ge j$, and with the convention that $\pmatrix{-1\cr -1\cr}=1$, 
$$t(i,j)=\sum_{f=0}^j\sum_{r=i-j}^{\min(f+1,i-f)}\pmatrix{f+1\cr r\cr}\pmatrix{i-f-1\cr r-1\cr}
\pmatrix{f+r+1\cr r-i+j\cr}\dnum$$

\entry{Remarks}  1.  One might hope to estimate $t(i,j)$ by using the naive Stirling approximations to the various binomial coefficients.  This turns out to be far too crude to be of any use.  

2.  The symbols $f,g,h,F,G$ and $H$ are hereby released from service and available to be used for other purposes in the remainder of the paper.

\newsec{The $t$-numbers by recursion.}

Let $f(i,j)$, $g(i,j)$ and $h(i,j)$ be the numbers of $(i,j)$-strings ending in [A], [B], and [AB].

An $(i,j)$ string ending in $[A]$ is constructed by appending an $[A]$ to
an $(i-1,j)$ string ending in $[B]$ or $[AB]$.  From this and similar
considerations we have

$$f(i,j)=g(i-1,j)+h(i-1,j)$$
$$g(i,j)=f(i,j-1)+h(i,j-1)$$
$$h(i,j)=f(i-1,j-1)+g(i-1,j-1)+h(i-1,j-1)$$
with initial conditions $f(1,0)=g(0,1)=h(1,1)=1$ and all other values
of all functions equal to zero whenever either $i$ or $j$ is non-positive.

Putting $F(x,y)=\sum f(i,j)x^iy^j$, $G(x,y)=\sum g(i,j) x^i y^j$,
$H(x,y)=\sum h(i,j)x^iy^j$, we solve to get

$$F(x,y)={x(1+y)\over 1-xy(2+x+y+xy)}$$
$$G(x,y)={(1+x)y\over 1-xy(2+x+y+xy)}$$
$$H(x,y)={x(1+x)y(1+y)\over 1-xy(2+x+y+xy)}$$

Clearly $t(i,j)=f(i,j)+g(i,j)+h(i,j)$, so if we set
$$T(x,y)=
\sum t(i,j)x^iy^j$$ we have $T=1+F+G+H$ and therefore
$$T(x,y)=   {1+x+y+xy\over 1-2xy-x^2y-xy^2-x^2y^2}\dnum$$

From this we see that $t$ is given by the following recursion relations:

\par
\centerline{\vbox{\vskip 8pt
\centerline{$$t(i,j)=0 \hbox{ if $i<0$ or $j<0$}$$}\vskip 2pt
\centerline{$$t(0,0)=t(0,1)=t(1,0)=1$$}\vskip 2pt
\centerline{$$t(i,0)=t(0,j)=0\hbox{ for $i,j\ge2$}$$}\vskip 2pt
\centerline{$$t(1,1)=t(1,2)=t(2,1)=3$$}\vskip 2pt
\centerline{$$\hbox{Otherwise }t(i,j)=t(i-2,j-2)+t(i-2,j-1)+t(i-1,j-2)+2t(i-1,j-1)$$}}}

In particular, we get the following table of values for $t(i,j)$:

$$\matrix{
&i=0&i=1&i=2&i=3&i=4&i=5&i=6&i=7&i=8&i=9&i=10\cr
j=0&1&1&0&0&0&0&0&0&0&0&0\cr
j=1&1&3&3&1&0&0&0&0&0&0&0\cr
j=2&0&3&9&10&5&1&0&0&0&0&0\cr
j=3&0&1&10&27&33&21&7&1&0&0&0\cr
j=4&0&0&5&33&83&108&81&36&9&1&0\cr
j=5&0&0&1&21&108&259&353&298&161&55&11\cr
j=6&0&0&0&7&81&353&817&1154&1066&665&281\cr
j=7&0&0&0&1&36&298&1154&2599&3776&3745&2612&\cr
j=8&0&0&0&0&9&161&1066&3776&8323&12371&12997\cr
j=9&0&0&0&0&1&55&665&3745&12371&26797&40586\cr
j=10&0&0&0&0&0&11&281&2612&12997&40586&86659\cr}$$

The main diagonal appears on the OEIS (up to a shift) 
as A171155, and the superdiagonal as A271943.

The generating function yields a crude upper bound on $t(i,j)$:

\entry{Theorem}
$$ t(i,j) \le \pmatrix{2j+1\cr i\cr}$$

{\bf Proof.}  Temporarily write $u(i,j)=\pmatrix{2j+1\cr i\cr}$.   Set
$$U(x,y)=
\sum_{i,j\ge 0}u(i,j)x^iy^j=\sum_{j\ge 0} y^j(1+x)^{2j+1}={1+x\over 1-y(1+x)^2}\dnum$$

Subtracting (\the\secno.1) from (\the\secno.\the\dnumno) gives
$$\sum_{i,j}\Big(u(i,j)-t(i,j)\Big) x^iy^j=U(x,y)-T(x,y)=
{y^2(1+x)^2\over \Big(1-y(1+x)^2\Big)\Big(1-xy(2+x+y+xy)\Big)}$$
and the power series of the right-hand side clearly has non-negative coefficients.  This proves the  inequality.

\qed

\newsec{Asymptotics}

In this section we will lay the foundations for studying the asymptotic behavior of $t(j+k,j)$ for fixed $k$ (that is, the northwest-to-southeast diagonals in the table above).

Recall that

$$T(x,y)=
\sum_{i,j} t(i,j)x^iy^j={(1+x)(1+y)\over 1-xy(2+x+y+xy)}\dnum$$

We want to estimate the $t(i,j)$.

\entry{Notation}   Set 
$$g(k,j)=t(k+j,j)\hbox{ or equivalently } g(i-j,j)=t(i,j)\dnum$$
and 
$$G_k(y)=\sum_{j=0}^\infty g(k,j)y^j \dnum$$

\entry{Proposition}
\itemitem{a)}$g(k,j)=g(-k,j+k)$
\itemitem{b)}$G_k(y)=G_{-k}(y)/y^k$
\itemitem{c)}$T(x,y/x)=\sum_{k=-\infty}^\infty x^k G_k(y)$

{\bf Proof.}  a) is an immediate consequence of the symmetry of $t$:
$$g(k,j)=t(k+j,j)=t(j,j+k)=g(-k,j+k)$$
b) is an immediate consequence of a) together with the displayed equation (\the\secno.2).

For c), write
$$\eqalign{T(x,y/x)&= \sum_{i,j} t(i,j)x^{i-j} y^j\cr
&=\sum_{i,j} t(k+j,j)x^k y^j \cr
&= \sum_{k} x^k G_k(y)\cr
}$$

\qed

\entry{Notation} Set:
$$Q= 1-4y+2y^2+y^4\qquad S=\sqrt{Q}\qquad A=1-2y-y^2\qquad M={A-S\over 2y^2}\dnum$$

\entry{Proposition}
\itemitem{a)}  $G_0(y)={1-y-S \over yS}$
\itemitem{b)}  For $k>0$, $G_{-k}(y)=y^kM^k\big(G_0(y)+1/y\big)$
\itemitem{c)}  For $k>0$, $G_k(y)=M^k\big(G_0(y)+1/y\big)$

{\bf Proof.}  For any fixed $k$, we have 
$$\eqalign{
{1\over 2\pi i}\int_{x\in C} {T(x,y/x)\over x^{1+k}}
&={1\over 2\pi i}\int_C\sum_m G_m(y) x^{m-k-1}\cr
&={1\over 2\pi i}\sum_m\left[G_m(y)\int_Cx^{m-k-1}\right]\cr
&=G_k(y)}$$
where $C$  is a small circle enclosing the origin.

Thus it suffices to compute the integral on the left for each $k$.
We can replace $y$ with a complex number close to the origin and restrict to an appropriately chosen annulus, so that the integral is a sum of residues at the  poles contained in $C$. From (\the\secno.1), one computes that (for any $k$) the poles are contained in the set
$$\left\{x=0,\quad x={A- S\over 2y},\quad x={A+S\over 2y}\right\}$$
By taking $y$ sufficiently small, we can ignore the pole with the plus sign, because it goes to $\infty$ as $y$ goes to $0$, and therefore can be taken to be outside the circle $C$.  

When $k=0$, this leaves the poles at $x=0$ and $x=(A-S)/(2y)$, where the residues are

$$R_1= -{1\over y}\quad \hbox{and}\quad R_2={(1-y)\over yS}$$
which proves a).

For $k<0$, the only pole is at $(A-S)/(2y)$ where the residue is $y^kM^k(1-y)/(yS)$, which, together with a), proves b).

For $k>0$, we get the same poles as in the case $k=0$ and can in principle proceed by adding the residues at both poles.  This, however, can get unpleasantly messy, so we proceed a different way by noting that  c) follows immediately from b) together with Proposition \the\secno.2(b).

\line{\hfill q.e.d.}

\newcount\msec\msec=\the\secno  \newcount\gsec\gsec=\the\secno 
\advance \msec by 2 \advance \gsec by 3

In Sections 5 and 6, we will approximate the power series coefficients of $G_k$, first for $k=0$ and then for $k\ge 1$.  First, we need some lemmas on approximating power series coefficients more generally.

\newsec{Some Approximation Lemmas}

\entry{Notation} The symbol $p$ always denotes an integer and $n$ always denotes a non-negative integer. For any such $p$ and $n$, set
$$c_p(n)=(-1)^n\pmatrix{p+{1\over2}\cr n\cr}$$

\entry{Lemma} 
There are positive constants $K_p$ and $L_p$ such that for all $n\ge 1$
 $$ L_p n^{-p-3/2}\le |c_p(n)|\le K_p{n^{-p-3/2}}$$

{\bf Proof.}   Note that
$$c_{p+1}(n)=\Bigg({p+3/2\over p+(3/2)-n}\Bigg)c_p(n)$$
and therefore it suffices to prove the lemma for a single value of $p$.   Taking $p=-1$, we have
$$c_{-1}(n)=(-1)^n\pmatrix{-1/2\cr n\cr}=4^{-n}\pmatrix{2n\cr n}$$
and
$${1\over \sqrt{(n+1)\pi}}\le 4^{-n}\pmatrix{2n\cr n}\le {1\over\sqrt{n\pi}}\dnum$$
as needed.

\qed

\entry{Remark}  The proof of Lemma \the\secno.2 explicitly gives upper and lower bounds for $\pmatrix{-1/2\cr n\cr}$.  It also implicitly gives the following upper and lower bounds:
$${1\over 2\sqrt{\pi}n^{3/2}}\le \left|\pmatrix{1/2\cr n\cr}\right| \le {1\over 2\sqrt{\pi}(n-1/2)^{3/2}}$$

\entry{Lemma} Let $\rho$ be a positive real number, $X=X(y)=1-y/\rho$, and $F(y)$ a function analytic on $|y|<r$ for some $r>|\rho|$, and $p$ an integer.  Set
$$\Phi(y)=\sum_n\phi(n)y^n=F(y)X^{p+1/2}$$

 Then
$$\phi(n)=F(\rho)c_p(n)\rho^{-n}+\bigo{\rho^{-n}n^{-p-5/2}}\dnum$$
with an implied constant that depends on $F$, $p$ and $r$.

{\bf Proof.}     Set $$V(y)={F(y)-F(\rho)\over y-\rho}=-{F(y)-F(\rho)\over \rho X}$$
Then $V$ is analytic on $|y|<r$.  Write $V(y)=\sum_j v_jy^j$.

 A little algebraic manipulation gives
$$\Phi(y)=F(\rho)X^{p+1/2}-\rho V(y)X^{p+3/2}$$
Taking coefficients, 
$$\phi(n)=F(\rho)c_p(n)\rho^{-n}-
\rho \sum_{j=0}^n  v_j c_{p+1}(n-j)\rho^{j-n}$$

Therefore (\the\secno.2) follows from
$$\sum_{j=0}^n  v_j c_{p+1}(n-j)\rho^{j-n}=\bigo{\rho^{-n}n^{-p-5/2}}$$
(with the fixed factor $\rho$ absorbed into the implied constant) or more simply
$$\sum_{j=0}^n  v_jc_{p+1}(n-j)\rho^{j}=\bigo{n^{-p-5/2}}\dnum$$

To prove (\the\secno.\the\dnumno), fix $R\in (\rho,r)$ and put $q=\rho/R<1$.   

$\underline{\hbox{\sl Claim 1.} }$  $v_j\rho^j=O(q^{j})$  Proof:  $V$ is analytic on $|y|<r$, so its Taylor series converges at $R<r$.  In particular its terms are bounded, so $v_j\rho^j=(v_jR_j)q^j=O(q^j)$.

$\underline{\hbox{\sl Claim 2.}}$  Uniformly for $0\le j\le n-1$,  $c_{p+1}(n-j)=\bigo{(n-j)^{-p-5/2}}$.  Proof:  This is Lemma \the\secno.2, with $p$ replaced by $p+1$.

$\underline{\hbox{\sl Claim 3.} }$ 
$\sum_{j=0}^{n/2} v_j c_{p+1}(n-j)\rho^j=\bigo{n^{-p-5/2}}$ (with an implied constant that might depend on $p$).  
Proof:  From Claims  1 and 2, the $j$th term is $v_j\rho^j c_{p+1}(n-j)=
\bigo{q^j}\bigo{(n-j)^{-p-5/2}}=
\bigo{q^j}\bigo{n^{-p-5/2}}$ 
with the last equality holding because $n-j>n/2$ throughout the relevant range.  Therefore (because $q<1$) the sum is 
$$\bigo{\Big(\sum_{j\ge 0} q^j\Big) n^{-p-5/2}}=\bigo{n^{-p-5/2}}$$

$\underline{\hbox{\sl Claim 4.}}$  $\sum_{j=n/2}^n
 v_j c_{p+1}(n-j)\rho^j=\bigo{n^{-p-5/2}}$.  
Proof:  In this range $j\ge n/2$ so $q^j\le q^{n/2}$.  
	Because $p$ is fixed, the   $c_{p+1}(n-j)$ are bounded in absolute value by a fixed polynomial in $n$.  Hence the entire sum is $\bigo{n^Aq^{n/2}}$ for some constant $A$.  Because exponential decay dominates every power of $n$, we get $n^Aq^{n/2}=\bigo{n^{-p-5/2}}$.

Claims 3 and 4 collectively imply (\the\secno.3) and therefore complete the proof.

\entry{Remark} The proofs of Lemmas \the\secno.2 and \the\secno.4 require no analytic machinery beyond elementary calculus.  In particular, the bounds in the displayed equation (\the\secno.1)
do not require Stirling's approximation.   They can be proved by elementary methods, and are so proved in Chapter 11 of [AN], via an appeal to Wallis's formula, which in turn is proved by elementary methods in [W].

Both lemmas are special cases of the transfer principle from singularity analysis (see Theorems VI.1 and VI.4 in [FS]), the invocation of which would be antithetical to the goals and the spirit of this paper.   The analyticity hypothesis in Lemma \the\secno.4 \hskip 2pt is stronger than the hypotheses of the [FS] theorems, which allows us to circumvent the hard analysis needed to establish the transfer principle in general.

\newsec{The Coefficients of $G_0$}

Now we will investigate the coefficients of $G_0=\sum_{n\ge 0}g(0,n)y^n$.

\entry{Notation} With $A$ and $Q$ as in Notation 3.4, let $\zeta$ be the root of $Q$ closest to the origin.  Readers with access to a good symbolic calculator will have no trouble verifying that 
$$\zeta\approx .296\qquad\hbox{and} \qquad\zeta^{-1}\approx 3.383$$

(Readers without access to a good symbolic calculator  can note that $Q=(1-y)(1-3y-y^2-y^3)$ and $\zeta$ is a root of the cubic factor, which is of course solvable in radicals.)  

Set
$$X=(1-y/\zeta)\qquad R=Q/X$$

and

$$C_1={1-\zeta\over 2\zeta^{2}\sqrt{2+\zeta}}\approx 2.660 \qquad C_2={C_1\over \sqrt{\pi}}\approx 1.501$$

All of this notation will remain in effect for the rest of the paper.

\entry{Lemma} $\sqrt{R}$ is analytic on an origin-centered disk of radius greater than $\zeta$.

{\bf Proof.}  Because $\zeta$ is both a simple root of $Q$ and the unique root closest to the origin, 
$R=Q/(1-y/\zeta)$ must be non-zero on such a disk.

\qed

\entry{Theorem}

\bigskip

\line{\hskip\parindent a)\hfill 
$\displaystyle g(0,n)=C_1\left|\pmatrix{-1/2\cr n+1\cr}\right|\times\zeta^{-n}+\bigo{\zeta^{-n}n^{-3/2}}$\hfill}
\bigskip
\line{\hskip \parindent b) \hfill $\displaystyle
g(0,n)\approx {C_2\times \zeta^{-n}\over \sqrt{n+1}}$\hfill}

{\bf Proof.}  Write 
$$B=\sum_{n\ge 0} b(n)y^{n}= yG_0+1 ={1-y\over \sqrt{Q}}={1-y\over\sqrt{R}}X^{-1/2}$$
Using Lemma 5.2, Lemma 4.4, and the easy calculation that $\sqrt{R(\zeta)}=2\zeta\sqrt{2+\zeta}$, conclude that
$$b(n)={1-\zeta\over 2\zeta\sqrt{2+\zeta}}\times\left|\pmatrix{-1/2\cr n\cr}\right|\times \zeta^{-n}+\bigo{\zeta^{-n}n^{-3/2}}$$
and therefore, for $n\ge 1$, 
$$g(0,n)={1-\zeta\over 2\zeta^2\sqrt{2+\zeta}}\times\left|\pmatrix{-1/2\cr n+1\cr}\right|\times \zeta^{-n}+\bigo{\zeta^{-n}n^{-3/2}}$$
as needed for a).

b) then follows from the approximation (4.1).

\qed

\newsec{The Coefficients of $G_k$ with $k\ge 1$}

\entry{Notation}  Continue to employ all the notation of 3.4 and 5.1.  Also set
$$P=2y^2 M=A-\sqrt{R}X^{1/2}$$

\entry{Lemma} $\displaystyle G_k=\left({1-y\over 2^k y^{2k+1}}\right)\times\left({1\over \sqrt{R}X^{1/2}}\right)\times P^k$

{\bf Proof.}  This follows from the definitions in Notation 3.4 and from Proposition 3.5a) and c).

\entry{Theorem}  Continue to use Notation 5.1.        
                                                                                Then for  fixed $k\ge 1$:                                                                           \bigskip                                                                                                                                                                                                                                        \line{\hskip\parindent a)\hfill
$\displaystyle g(k,n)=C_1\left|\pmatrix{-1/2\cr n+2k+1\cr}\right|\times\zeta^{-n-k/2}+\bigo{\zeta^{-n}n^{-3/2}}$\hfill} 
with an implied constant that depends on $k$.
                                               \bigskip
\line{\hskip \parindent b) \hfill $\displaystyle
g(k,n)\approx {C_2\times \zeta^{-n-k/2}\over \sqrt{n+2k+1}}$\hfill}

{\bf Proof.}  Start with the binomial expansion
$$P^k=\sum_{j=0}^k  (-1)^j\pmatrix{k\cr j\cr} A^{k-j} \sqrt{R}^jX^{j/2}$$

Then from Lemma \the\secno.2, we have
$$2^k y^{2k+1}G_k=\sum_{j=0}^k (1-y)(-1)^j\pmatrix{k\cr j\cr}A^{k-j}\sqrt{R}^{j-1}X^{(j-1)/2}\dnum$$

Consider the sum on the right side of (\the\secno.\the\dnumno).

Then for fixed $k$, we have:

$\underline{\hbox{\sl Observation 1:\/}}$  The terms 
with odd $j$  have a finite number of non-zero coefficients.

$\underline{\hbox{\sl Observation 2:\/}}$   By Lemmas 5.2 and 4.4, the  terms  with even $j>0$  are all of the form $X^{1/2}\times\hbox{(analytic function)}$ and hence have coefficients that are $O\Big(\zeta^{-n}n^{-3/2}\Big)$

$\underline{\hbox{\sl Observation 3:\/}}$  By Proposition \the\secno.2 and Lemma 4.4, the coefficients of the $j=0$ term are
$$
(1-\zeta)A(\zeta)^{k}
\left(\sqrt{R(\zeta)}\right)^{-1}c_{-1}(n)\zeta^{-n}+\bigo{\zeta^{-n}n^{-3/2}}\dnum$$

Together, these observations  imply that the coefficients of $(\the\secno.1)$ are also of the form \the\secno.\the\dnumno.

The coefficients of $2^kG_k$ are then given by (\the\secno.2) with an offset of $2k+1$:
$$
\eqalignno{2^kg(k,n)=
(1-\zeta)A(\zeta)^{k}
\left(-\sqrt{R(\zeta)}\right)^{-1}c_{-1}(n+2k+1)&\zeta^{-n-2k-1}\cr
&+\bigo{\zeta^{-n-2k-1}n^{-3/2}}&(6.3)\cr}$$

One easily computes that $A(\zeta)=2\zeta^{3/2}$ and $\sqrt{R(\zeta)}=2\zeta\sqrt{2+\zeta}$, 
so (\the\secno.\the\dnumno) can be rewritten to establish a) and b) then follows from a) together with the approximation (4.1).

\qed

\entry{Remark} The statement of  \the\secno.3, when specialized to $k=0$, is identical to the statement of Theorem 5.3.  However, the proof of Theorem 6.3 does not work when $k=0$ (because Proposition 3.5c) does not hold for $k=0$), which is why we have two separate theorems.  (The proofs do, however, have several elements in common.)

\newsec{The Coefficients $t(i,j)$}

Recall from Notation 3.1 that for $i\ge j$, $t(i,j)=g(i-j,j)$.

\entry{Theorem}  Let $C_1$ and $C_2$ be as in Notation 5.1.  Then for $i\ge j\ge 1$  and for fixed $k=i-j$, we have

\bigskip                                                                                                                                                                                                                                        \line{\hskip\parindent a)\hfill
$\displaystyle t(i,j)=C_1\left|\pmatrix{-1/2\cr 2i-j+1\cr}\right|\times\zeta^{-(i+j)/2}+\bigo{\zeta^{-j}j^{-3/2}}$\hfill}                                           
with an implied constant that depends on $k$.     \bigskip
\line{\hskip \parindent b) \hfill $\displaystyle
t(i,j)\approx {C_2\times \zeta^{-(i+j)/2}\over \sqrt{2i-j+1}}$\hfill}

{\bf Proof.} This is a straightforward translation of Theorems 5.3 and 6.3.

\qed

We will also need a uniform bound on $t(i,j)$:

\entry{Theorem} $\displaystyle t(i,j)=O(\zeta^{-(i+j)/2})\hbox{ uniformly in $i$ and $j$ }$

{\bf Proof.}  Put $s(p)=\sum_{i+j=p}t(i,j)$ and define
$$S(z)=\sum_{p=0}^\infty s(p)z^p=T(z,z)={(1+z)^2\over 1-2z^2-2z^3-z^4}$$
Check that $\sqrt{\zeta}$ is a root of the denominator of $S$, and indeed the root of smallest modulus, so partial fractions give
$$s(p)=O\Big(\zeta^{-p/2}\Big)$$
Now if we set $p=i+j$, then $t(i,j)\le s(p)$ so the result follows.

\qed

\bigskip\bigskip

\centerline{\bf PART TWO:  HOTELLING PROTOCOLS}

We can now return to the question we started with.

Alice and Bob play 
 a game.  Alice has $m$ strategy sets and Bob has $n$.  We say that two moves from the same strategy set have the same {\it type\/}.  The example to keep in mind is that the first strategy set consists of prices, the second consists of locations, and so on.

Over the course of the game, Alice makes exactly one move from each of her strategy sets and likewise for Bob.  Payoffs depend on the strategies they play, and not on the order in which they play them.  

A {\it Hotelling protocol\/} is a list of instructions to the players that tells them which strategy sets to play from at each turn.  An example of a Hotelling protocol is ``At Turn 1, Alice announces location.  At Turn 2, Alice announces price and Bob announces location.  At Turn 3, Bob announces price.''  

We require that at each turn, at least one player makes at least one move.  We require also that there are not two turns in a row on which Alice (or Bob) does not play.  Two such turns, after all, might as well be combined into a single turn.

Our goal here is to estimate the number of Hotelling protocols as a function of $m$ and $n$.  

\newsec{The $u$-numbers.}

A Hotelling protocol (determining which moves each player makes at each turn) consists of a simple protocol (determining who plays at each turn) together with a specification of move types for each player at each turn.  

We continue to suppose that Alice has $m$ moves and Bob has $n$.  

To construct a Hotelling protocol, we can proceed as follows:

1)  Construct a simple protocol for a game in which Alice has $i$ moves and Bob has $j$ moves,
with $i\le m$ and $j\le n$.   Remember that a simple protocol can be represented as a sequence of symbols, each of which is either $[A]$, $[B]$ or $[AB]$, as in Definition 1.1.   Call each $[A]$ and each $[AB]$ an {\sl occurrence\/} of $A$, and likewise for $B$.

2)  Assign each of Alice's $m$ moves to one of the $i$ occurrences of $A$ and each of Bob's $n$ moves to one of the $j$ occurrences of $B$, with each occurrence receiving at least one move.

The number of ways to do this is given by
$$i!j!{m\brace i}{n\brace j}$$
where the braces indicate Stirling numbers of the second kind.

Therefore we have 

$$u(m,n)=\sum_{\vbox{\hbox{$i\le m$}\hbox{$j\le n$}}}
i!j!{m\brace i}{n\brace j}t(i,j)\dnum$$

where $t(i,j)$ is as defined in Definition 1.2 and approximated in Theorem 7.1.

These numbers grow quite fast, as indicated by the following table:

{
\font\seven=cmr7
\textfont0=\seven  \fam0

$${\fam0
\matrix{&m=1&m=2&m=3&m=4&m=5&m=6&m=7&m=8\cr n= 1  &
3& 9& 27& 81& 243& 729& 2187& 6561\cr n=2   & 
9& 51& 261& 1299& 6429& 31851&
   158181& 787299\cr n= 3  & 
27& 261& 2067& 15165& 107547& 751221& 5215587& 
  36152685\cr n= 4  & 
81& 1299& 15165& 155187& 1487541& 13769259& 125038365& 
  1123781187\cr n= 5  & 
243& 6429& 107547& 1487541& 18648483& 221003949& 
  2531982507& 28409691621\cr n= 6  & 729& 31851& 751221& 13769259& 221003949& 
  3279778851& 46332606741& 633837539259\cr n=7   & 
2187& 158181& 5215587& 
  125038365& 2531982507& 46332606741& 794322896307& 
  13033989543885\cr n= 8  & 
6561& 787299& 36152685& 1123781187& 28409691621& 
  633837539259& 13033989543885& 253477135842387\cr}}$$
  
 }

\newsec{Estimating $u(m,n)$}

To estimate $u(m,n)$, we will allow ourselves to call on deeper analytical tools than we allowed in Part One.  Specifically, we need Temme's approximations for Stirling numbers, in the form of the facts listed in \the\secno.4 below.  First we need to define some constants and some additional notation.

\entry{Constants} 
As in earlier sections, $\zeta\approx .295598$ is the real root of  $y^3+y^2+3y-1$.  Also set:

$$r_0=\log\Big(1+\sqrt{\zeta}\Big)\approx .434175\qquad
\alpha_0={\sqrt{\zeta}\over r_0 (1+\sqrt{\zeta})}\approx .811196\qquad Z=1/r_0\approx 2.30322$$

\entry{Notation}
$$\alpha\mapsto r(\alpha) \hbox{ is the function defined implicitly by the equation } \alpha={1-e^{-r}\over r}\dnum$$

$$f(\alpha)=\alpha\log\Big(e^{r(\alpha)}-1\Big)-\log(r(\alpha))\dnum$$

$$\widehat{t}(i,j)=\zeta^{-(i+j)/2}\approx 1.83929^{i+j}$$

$$\Psi_{m}(i)=i!{m\brace i}\zeta^{-i/2}$$

$$\Theta_m=\sum_{i\le m}\Psi_m(i)$$

$$\Phi_{m}(i)=mf(i/m)+i \log(\zeta^{-1/2})$$

$$\widehat{u}(m,n)=
\sum_{\vbox{\hbox{$i\le m$}\hbox{$j\le n$}}}
i!j!{m\brace i}{n\brace j}\widehat{t}(i,j)=\Theta_m\times \Theta_n$$

\entry{Definition}  For sequences $(a(n))$ and $(b(n))$, I will write
$$a(n)=b(n) n^{O(1)}$$
to mean that there is a constant $C>0$ such that
$$n^{-C}\le {a_n\over b_n}\le n^C\dnum$$
For doubly indexed sequences $(a(i,n))$ and $(b(i,n))$ I will say that
$a(i,n)=b(i,n) n^{O(1)}$  uniformly in $i$
if there is a single constant $C$ with $n^{-C}\le a(i,n)/b(i,n)\le n^C$ for all $i$.

\entry{Some Facts} 

a) $\displaystyle f'(\alpha)=\log\Big(e^{r(\alpha)}-1\Big)$

b)  $\displaystyle i! {n\brace i}= n!e^{n f(i/n)}n^{O(1)}$ uniformly for $i\le n$.

{\bf Proof.}  See [T], sections 2 and 3.

\entry{Proposition} $u(m,n)=O\Big(\widehat{u}(m,n)\Big)$

{\bf  Proof.}  This follows immediately from the definitions and Theorem 7.2.

\qed

In view of the proposition, our goal is to estimate $\widehat{u}(m,n)=\Theta_m\times \Theta_n$.  For this it will suffice to estimate $\Theta_m$, and for this we will start by estimating $\Psi_m(i)$.

\entry{Proposition}
$\displaystyle\log\Big(\Psi_m(i)\Big)=\log(m!)+\Phi_m(i)+O(\log(m))$

{\bf Proof.}  This follows immediately from Fact \the\secno.4b).

\entry{Theorem}  

a)  As $m\rightarrow\infty$, 
$$\max_{1\le i \le m}\Psi_m(i)=m!Z^m m^{O(1)}$$

b)  The exponential part of $\Psi_m(i)$ is uniquely maximized at $i\approx \alpha_0m$.

{\bf Proof.}  From Proposition \the\secno.6, it suffices to maximize $\Phi_{m}(i)=mf(i/m)+i\log(1/\sqrt{\zeta})$.   Setting $i/m=\alpha$, this becomes
$$mf(\alpha)+m\alpha\log(1/\sqrt{\zeta})\dnum
$$
The first-order condition for an optimum is
$$f'(\alpha)=-\log(1/\sqrt{\zeta})\dnum$$

Now Fact \the\secno.4a) allows us to conclude that if $\alpha$ is chosen optimally,   $r(\alpha)=r_0$ and the displayed equation (\the\secno.1) allows us to conclude that the optimum occurs at $\alpha_0$, which establishes b).

Setting $\alpha=\alpha_0$ in the maximand (\the\secno.3) and using the defining equation for $f$ (\the\secno.2) we find that at the optimum $i_0=m\alpha_0$ we have 
$$
\Phi_{m}(i_0)=
mf(\alpha_0)+\alpha_0 m\log(\sqrt{1/\zeta})=m\log(Z)$$
which, together with Proposition \the\secno.6, establishes a).

\qed

\entry{Corollary}

 a) $\displaystyle \widehat{u}(n,n)= (n!)^2\times Z^{2n} \times n^{O(1)}\approx (n!)^2\times 5.30482^n\times n^{O(1)}$

b)   $\displaystyle {u}(n,n)= (n!)^2\times Z^{2n} \times n^{O(1)}\approx (n!)^2\times 5.30482^n\times n^{O(1)}$

{\bf Proof.} a)  $\widehat{u}(n,n)$ is a sum of $n^2$ terms, so (because $n^2$ is a polynomial) its exponential growth is determined by the exponential growth of the fastest-growing term.   Theorem \the\secno.7 completes the proof of a).

b)  From part a) we have $u(n,n)\le (n!)^2 Z^{2n} n^{O(1)}$.   For the opposite inequality, let $i_0$ be the integer closest to $\alpha_0n$ (so that $i_0$ depends on $n$) and make the following observations:

$\underline{\hbox{\sl Observation 1:\/}}$ By Theorem \the\secno.7, together with the sub-observation that moving from $\alpha_0m$ to the nearest integer $i_0$ changes $\Phi_n$ by an amount that can be easily absorbed into $n^{O(1)}$, we have: 
$$\Psi_n(i_0)=i_0!{n \brace i_0}\zeta^{-i_0/2}\ge n! Z^n n^{O(1)}$$

$\underline{\hbox{\sl Observation 2:\/}}$ 
By Theorem 7.1,
$$t(i_0,i_0)=\zeta^{-i_0}n^{O(1)}$$

$\underline{\hbox{\sl Observation 3:\/}}$  Because all terms in (8.1) are non-negative,
$$u(n,n)\ge\left(
i_0!{n \brace i_0}\right)^2 t(i_0,i_0)$$

$\underline{\hbox{\sl Observation 4: \/}}$ Combining the three preceding observations, we have 
$$u(n,n)\ge (n!)^2 Z^{2n} n^{O(1)}$$

Combined with Proposition \the\secno.5  and part a), this proves part b).

\qed

\noindent{\bf Remarks on the case $m\neq n$:}
If $m\neq n$, Proposition 9.5 still gives
$$u(m,n)\le m!n! Z^{m+n} (m+n)^{O(1)}$$

To obtain a matching lower bound by arguments analogous to those of Corollary \the\secno.8, we would need an estimate for $t(i,j)$ that remains valid up to polynomial factors when $i-j$ grows proportionally with $i+j$.  Theorem 7.1 does not provide this because its implied constant is allowed to vary with $k=i-j$.  Therefore it does not control the terms near $i\approx \alpha_0m, j\approx \alpha_0n$ when $m-n$ itself grows.  I don't know whether the proposed matching lower bound nevertheless holds.

\newsec{Concluding Remarks}

Set 
$$\eqalign{H(n)&=u(n,n)\cr
&=\hbox{ The number of Hotelling protocols when Alice and Bob have $n$ moves each}\cr
\widehat{H}(n)&=(n!)^2\times 5.30482^n\cr} $$

Corollary 9.7b) establishes that $\widehat{H}(n)$ is a good approximation to $u(n,n)$ in the sense that
$$H(n)=\widehat{H}(n)\times n^{O(1)}$$

The following graph illustrates
the values of $\log(H(n))$ (in blue) for even values of $n$ and of $\log(\widehat{H}(n))$ (in gold) for odd values of $n$  in the range $1\le n\le 100$.   Between $n=1$ and $n=100$, the difference 
$\log(\widehat{H}(n))-\log(H(n))$ ranges from .570 to 2.623.

\bigskip

\centerline{\includegraphics[width=5in]{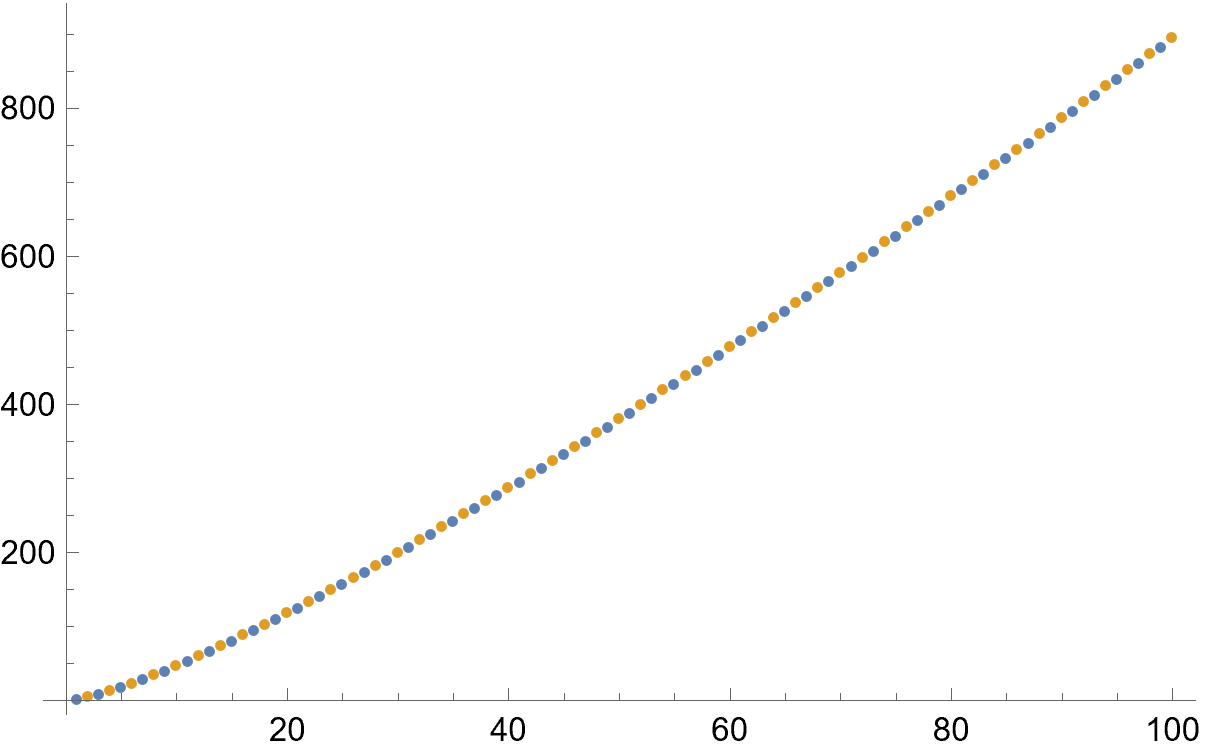}}
 
\vfill\eject

\centerline{\bf REFERENCES }

[AN]  Alsina, C., \& Nelsen, R. B. (2010). {\it Charming Proofs: A Journey 
into Elegant Mathematics},  Mathematical Association of America.

[FS]  Flajolet, P. and Sedgewick, R., {\it Analytic Combinatorics\/}, Cambridge University Press (2009)

[MO]  MathOverflow,  ``Approximating power series coefficients --- Why does a clearly illegitimate method (sometimes) work so well?", https://mathoverflow.net/questions/364319

[PW]  Pemantle, R,  Wilson, M., and Melczer, S., {\it Analytic Combinatorics in Several Variables\/} (2nd edition), Cambridge University Press 2024.

[T]  Temme, N., ``Asymptotic estimates of Stirling numbers'',
  {Studies in Applied Mathematics} 89, (1993), 233-243.

[W] ``An Elementary Proof of the Wallis Product Formula for $\pi$" by Johan Wastlund, published in The American Mathematical Monthly, Vol. 114, No. 10 (December 2007), pp. 914-917
\bye